%% file: agt-1-34.tex
\documentclass{gtart}
\input agtout

\lognumber{34}
\volumenumber{1}
\volumeyear{2001}
\papernumber{34}
\pagenumbers{699}{708}
\received{2 August 2001}
\revised{15 November 2001}
\accepted{16 November 2001}
\published{22 November 2001}

\usepackage{amsmath,amssymb}
\usepackage{curvesls}
\usepackage{epic}

\newtheorem{thm}{Theorem}[section]
\newtheorem{lem}[thm]{Lemma}

\newtheorem{prop}[thm]{Proposition}

\theoremstyle{definition}

\newcommand{\C}{{\bf{C}}}

\newcommand{\Zed}{{\bf{Z}}}
\newcommand{\GL}{{\mathrm{GL}}}
\newcommand{\Diff}{\mathrm{Diff}}
\newcommand{\DiffZ}{{\mathrm{Diff}^{\Zed_2}}}
\newcommand{\Map}{{\mathrm{Map}}}

\newcommand{\GTS}{{\Sigma_2}}

\newcommand{\Stab}[1]{\mathrm{Stab}({#1})}
\newcommand{\ab}{{\mathrm{ab}}}

\renewcommand{\exp}{{\mathrm{exp}}}

\begin{document}

\title{The mapping class group of a genus two\\surface is linear}
\asciititle{The mapping class group of a genus two surface is linear}
\author{Stephen J. Bigelow\\Ryan D. Budney}
\asciiaddress{Department of Mathematics and Statistics, University of 
Melbourne\\Parkville, Victoria, 3010, Australia\\and\\Department 
of Mathematics, Cornell University\\Ithaca, 
New York 14853-4201, USA}
\address{Department of Mathematics and Statistics, University of 
Melbourne\\Parkville, Victoria, 3010, Australia}
\secondaddress{Department 
of Mathematics, Cornell University\\Ithaca, 
New York 14853-4201, USA}
\email{bigelow@unimelb.edu.au}
\secondemail{rybu@math.cornell.edu}
\asciiemail{bigelow@unimelb.edu.au, rybu@math.cornell.edu}
\begin{abstract}
In this paper we construct a faithful representation of the mapping
class group of the genus two surface into a group of matrices over the
complex numbers.  Our starting point is the Lawrence-Krammer
representation of the braid group $B_n$, which was shown to be
faithful by Bigelow and Krammer.  We obtain a faithful representation
of the mapping class group of the $n$-punctured sphere by using the
close relationship between this group and $B_{n-1}$.  We then extend
this to a faithful representation of the mapping class group of the
genus two surface, using Birman and Hilden's result that this group is
a $\Zed_2$ central extension of the mapping class group of the
$6$-punctured sphere.  The resulting representation has dimension
sixty-four and will be described explicitly. In closing we will remark
on subgroups of mapping class groups which can be shown to be linear
using similar techniques.
\end{abstract}

\asciiabstract{In this paper we construct a faithful representation of
the mapping class group of the genus two surface into a group of
matrices over the complex numbers.  Our starting point is the
Lawrence-Krammer representation of the braid group B_n, which was
shown to be faithful by Bigelow and Krammer.  We obtain a faithful
representation of the mapping class group of the n-punctured sphere
by using the close relationship between this group and B_{n-1}.  We
then extend this to a faithful representation of the mapping class
group of the genus two surface, using Birman and Hilden's result that
this group is a Z_2 central extension of the mapping class group
of the 6-punctured sphere.  The resulting representation has
dimension sixty-four and will be described explicitly. In closing we
will remark on subgroups of mapping class groups which can be shown to
be linear using similar techniques.}

\primaryclass{20F36}
\secondaryclass{57M07, 20C15} 
\keywords{Mapping class group, braid group, linear, representation} 

\maketitle


\section{Introduction}
\label{INTRODUCTION}

Let $\Diff M$ denote the topological group of 
orientation preserving diffeomorphisms 
of an oriented manifold $M$
which act as the identity on $\partial M$.
The {\em mapping class group} of $M$
is the group $\pi_0 \Diff M$.
A {\em representation} of a group
is a homomorphism from the group
into a multiplicative group of matrices over some commutative ring.
A representation is called {\em faithful} if it is one-to-one.
A group is called {\em linear} if it admits a faithful representation.

The aim of this paper is to construct
a faithful representation of the mapping class group
of the genus two surface. In the process we construct
faithful representations of mapping class groups of 
punctured spheres, hyperelliptic mapping class groups 
and, more generally, normalizers of certain 
covering transformation groups of surfaces. 

We take as our starting point
the Lawrence-Krammer representation
of the braid group $B_n$.
Bigelow \cite{Bigelow} and Krammer \cite{Krammer2}
have shown this to be faithful.
In Section \ref{SPHERE},
we show how to alter the Lawrence-Krammer representation
to obtain a faithful representation
of the mapping class group of an $n$-times punctured sphere.

The genus two surface is a branched covering space of the sphere
with six branch points.
Birman and Hilden \cite{Birman}
have used this fact to establish a close relationship
between the mapping class group of the genus two surface
and the mapping class group of the six-times punctured sphere.
In Section \ref{GENUS2},
we use this relationship to obtain a faithful representation
of the mapping class group of the genus two surface.

Simultaneous with this result, Nathan Dunfield and also Mustafa Korkmaz 
\cite{Korkmaz}
have individually produced faithful
representations of the mapping class group of the genus two surface.
All of these constructions use the relationship to
the mapping class group of the six-times punctured sphere.
However we have taken a bit of extra care to keep the dimension reasonably low.
Our faithful representations of
the mapping class groups of
the $n$-times punctured sphere and the genus two surface
have dimensions $n \binom{n-1}{2}$ and $64$ respectively,
whereas the representations in \cite{Korkmaz} have dimensions
$n \binom{n-1}{2}^2$ and $2^{10}3^5 5^3$ respectively.

The low rank of our
representation makes it suitable for computer use, and we explicitly compute
the matrices for our representations in Section \ref{MATRICES}.
In Section \ref{REMARKS}
we show how to generalize our construction to obtain 
faithful representations normalizers of a class of 
finite subgroups of mapping class groups.
The simplest such generalization
gives a faithful representation
of the hyperelliptic group of the genus $g$ surface.
Korkmaz \cite{Korkmaz} also constructed a faithful representation
of the hyperelliptic group,
but once again ours has a smaller dimension,
namely $(2g+2)\binom{2g+1}{2}+2g$
as opposed to $(2g+2)\binom{2g+1}{2}^2 3^{g^2} \prod_{i=1}^g(3^{2i}-1)$.

Throughout this paper,
$D$ will denote a disk,
$\GTS$ will denote a closed oriented surface of genus two,
and $S^2$ will denote a sphere.
If $M$ is an oriented manifold
and $n$ is a positive integer then
let $\Diff(M,n)$ denote $\Diff(M, \{p_1,\dots,p_n\})$,
where $p_1,\dots,p_n$ are distinct points in the interior of $M$.
This is the group of diffeomorphisms of $M$ 
that restrict to permutations of the set $\{p_1, \dots, p_n\}$.

\newpage
\section{The $n$-punctured Sphere}
\label{SPHERE}

The aim of this section is to prove the following.

\begin{thm}
\label{sphere}
There exists a faithful representation of 
the mapping class\break group of the $n$-times punctured sphere.
\end{thm}

The braid group $B_n$ is the group $\pi_0 \Diff(D,n)$.
Provided $n \ge 3$,
the center of $B_n$ is isomorphic to $\Zed$
and is generated by the {\em full twist braid} $\Delta^2$.
This is a Dehn twist about a curve which is parallel to $\partial D$.

Let $p_1,\dots,p_n$ be distinct points in $S^2$.

\begin{lem}
\label{disk_sphere}
Provided $n \ge 4$, there is a short exact sequence
\[
0\rightarrow \Zed
\rightarrow B_{n-1}
\rightarrow \Stab{p_n}
\rightarrow 0,
\]
where the image of $\Zed$ in $B_{n-1}$ is the center of $B_{n-1}$,
and
$\Stab{p_n}$ is the subgroup of $\pi_0 \Diff(S^2,n)$ consisting of 
diffeomorphisms that fix the point $p_n$.
\end{lem}

\begin{proof}
Let $D^+$ and $D^-$ be
the northern and southern hemispheres of $S^2$, that is, 
two disks in $S^2$ such that $D^+ \cap D^- = \partial D^+ = \partial D^-$.
Assume that $p_1,\dots,p_{n-1} \in D^+$ and $p_n \in D^-$.
Then $B_{n-1}$ is $\pi_0 \Diff(D^+,n-1)$.
We can extend any $f \in \Diff(D^+,n-1)$
to a diffeomorphism of the whole sphere
by setting it to be the identity on $D^-$.
Let $\phi \co B_{n-1} \rightarrow \pi_0 \Diff(S^2,n)$ be the homomorphism
defined in this way. 
This will be the rightmost map in our short exact sequence.

First we show that the image of $\phi$ is $\Stab{p_n}$.
Let $g$ be an element of $\Diff(S^2,n)$ which fixes the puncture $p_n$.
Note that $g_{|(D^-)}$ is a closed tubular neighborhood of $p_n$ in 
$S^2-\{p_1,\cdots,p_{n-1}\}$. By the uniqueness of tubular neighborhoods
theorem, $g_{|(D^-)}$ is isotopic to the identity relative to $\{p_n\}$. 
This isotopy can be extended to an ambient isotopy of 
the $n$-times punctured sphere.
We can therefore assume, without loss of generality,
that $g$ acts as the identity on $\partial D^-$.
Thus $g = \phi(g_{|(D^-)})$.

Now we show that the kernel of $\phi$ is generated by $\Delta^2$.
Let $f \in \Diff(D^+,n-1)$ represent an element of the kernel of $\phi$.
Let $g = \phi(f)$ be its extension to $S^2$ which is the identity on
$D^-$. Then there is an isotopy $g_t\in \Diff(S^2,n)$ such that
$g_0 = g$ and $g_1$ is the identity map. 
Now $g_t$ restricted to $D^-$ defines an element of the fundamental group of
the space of all tubular neighborhoods of $p_n$.
The proof of the uniqueness of tubular neighborhoods theorem 
\cite{Hirsch} naturally extends to a proof that there is a homotopy equivalence
between the space of tubular neighbourhoods of a point and $\GL(T_{p_n})$.
Thus the fundamental group of the space of tubular neighbourhoods of a 
fixed point in $S^2$ is $\Zed$,
generated by a rigid rotation through an angle of $2\pi$.
Consequently our family of diffeomorphisms $g_t$
can be isotoped relative to endpoints
so that its restriction to $D^-$
is a rigid rotations by some multiple of $2\pi$.
Therefore $f$ is isotopic to some power of $\Delta^2$.
\end{proof}

Let\vglue -40pt
\[
\mathcal L_n \co B_n \rightarrow 
\GL(\binom{n}{2},\Zed[q^{\pm 1},t^{\pm 1}])
\]
denote the Lawrence-Krammer representation,
which was shown to be faithful in \cite{Bigelow} and \cite{Krammer2}.
By assigning algebraically independent complex values to $q$ and $t$,
we consider the image as lying in $\GL(\binom{n}{2},\C)$.

Now $\mathcal L_n(\Delta^2)$ is a scalar matrix $\lambda I$.
This can be seen by looking at
the representation as an action on the module of forks \cite{Krammer1}.
(In fact, $\lambda = q^{2n}t^2$.)
We will now ``rescale'' the representation $\mathcal L_n$
so that $\Delta^2$ is mapped to the identity matrix.

The abelianization of $B_n$ is $\Zed$.
Let $\ab \co B_n \rightarrow \Zed$
denote the abelianization map.
Then $\ab(\Delta^2) \neq 0$,
as is easily verified using the standard group presentation for $B_n$.
(In fact, $\ab(\Delta^2) = n(n-1)$.)
Let $\exp \co \Zed \rightarrow \C^*$
be a group homomorphism which takes $\ab(\Delta^2)$ to $\lambda^{-1}$.
We now define a new representation $\mathcal L_n'$ of $B_n$ by
$$\mathcal L_n'(\beta) = (\exp \circ \ab (\beta))\mathcal L_n(\beta).$$
We claim that the kernel of $\mathcal L_n'$
is precisely the center of $B_n$, provided $n \ge 3$.
By design,
$\mathcal L_n'(\Delta^2) = I$.
Conversely, suppose $\mathcal L_n'(\beta) = I$.
Then $\mathcal L_n(\beta)$ is a scalar matrix, so lies in
the center of the matrix group.
Since $\mathcal L_n$ is faithful, it follows that
$\beta$ lies in the center of the braid group.

We are now ready to prove Theorem \ref{sphere}.
If $n \le 3$ then $\Diff(S^2,n)$
is simply the full symmetric group on the puncture points,
so the result is trivial.
We therefore assume $n \ge 4$.
By Lemma \ref{disk_sphere},
$\mathcal L_{n-1}'$ induces a faithful representation of $\Stab{p_n}$.
Since $\Stab{p_n}$ has finite index in $\pi_0 \Diff(S^2,n)$,
$\mathcal L_{n-1}'$ can be extended to a finite dimensional representation
$\mathcal K_n$ of $\pi_0 \Diff(S^2,n)$.
Extensions of faithful representations are faithful
(see for example \cite{Lang}), giving the result.

Note that the faithful representation $\mathcal K_n$
has dimension $n\binom{n-1}{2}$.

\eject
\section{The genus two surface}
\label{GENUS2}

The aim of this section is to prove the following.

\begin{thm}
There exists a faithful representation of 
the mapping class\break group of the genus two surface.
\end{thm}

\begin{figure}[ht!]
\centering

\begin{picture}(280,40)(0,0)

\thinlines

\qbezier(20,20)(20,40)(50,40)
\qbezier(50,40)(54,40)(60,37)
\qbezier(60,37)(70,32)(80,37)
\qbezier(80,37)(86,40)(90,40)
\qbezier(90,40)(120,40)(120,20)

\qbezier(20,20)(20, 0)(50, 0)
\qbezier(50, 0)(54, 0)(60, 3)
\qbezier(60, 3)(70, 8)(80, 3)
\qbezier(80, 3)(86, 0)(90, 0)
\qbezier(90, 0)(120, 0)(120,20)

\qbezier(40,20)(50,10)(60,20)
\qbezier(40,20)(50,30)(60,20)
\qbezier(40,20)(39,21)(38,22)
\qbezier(60,20)(61,21)(62,22)

\qbezier(80,20)(90,10)(100,20)
\qbezier(80,20)(90,30)(100,20)
\qbezier(80,20)(79,21)(78,22)
\qbezier(100,20)(101,21)(102,22)

\qbezier(180,20)(180, 0)(210, 0)
\qbezier(210, 0)(214, 0)(220, 3)
\qbezier(220, 3)(230, 8)(240, 3)
\qbezier(240, 3)(246, 0)(250, 0)
\qbezier(250, 0)(280, 0)(280,20)

\qbezier(180,20)(180,25)(190,25)
\qbezier(190,25)(195,25)(200,20)
\qbezier(200,20)(210,10)(220,20)
\qbezier(220,20)(230,30)(240,20)
\qbezier(240,20)(250,10)(260,20)
\qbezier(260,20)(265,25)(270,25)
\qbezier(270,25)(280,25)(280,20)

\put(180,20){\circle*{4}}
\put(200,20){\circle*{4}}
\put(220,20){\circle*{4}}
\put(240,20){\circle*{4}}
\put(260,20){\circle*{4}}
\put(280,20){\circle*{4}}

\multiput(0,20)(40,0){3}{\line(1,0){20}}
\multiput(20,20)(40,0){3}{\dashbox{1}(20,0){}}
\put(120,20){\line(1,0){14}}
\put(136,20){\line(1,0){4}}

\qbezier(132.89,15)(135,17.32)(135,20)
\qbezier(135,20)(135,30)(130,30)
\qbezier(130,30)(127.11,30)(125.67,25)
\put(126,25){\vector(-1,-3){0}}

\put(150,20){\line(1,0){20}}
\drawline(170,20)(167,23)
\drawline(170,20)(167,17)
\end{picture}
\caption{The action of $\Zed_2$ on $\GTS$.}
\label{oneeighty}
\end{figure}
The {\em standard involution} of $\GTS$
is the rotation through an angle of $\pi$
as shown in Figure \ref{oneeighty}.
This defines an action of $\Zed_2$
as a group of branched covering transformations
with quotient $S^2$ and six branch points.
Let $\DiffZ \GTS$ denote
the group of $\Zed_2$-equivariant diffeomorphisms of $\GTS$,
that is, the group of diffeomorphisms which
strictly commute with the standard involution.
We think of $\DiffZ\GTS$ as a subspace of $\Diff\GTS$.

\begin{prop}
The inclusion map
$\DiffZ\GTS \rightarrow \Diff\GTS$
induces an isomorphism on $\pi_0$.
\end{prop}

\begin{proof}
That the induced map is epic follows from Lickorish's theorem \cite{Lickr} 
that that the genus two mapping class group is generated by 
five Dehn twists, all of which happen to be $\Zed_2$ equivariant. 
See Figure \ref{genericlickgen}.  This is
the point where the analogous theorem fails for higher genus surfaces.
\begin{figure}[ht!]
\centering
\begin{picture}(100,40)(20,0)

\thicklines

\qbezier(20,20)(20,40)(50,40)
\qbezier(50,40)(54,40)(60,37)
\qbezier(60,37)(70,32)(80,37)
\qbezier(80,37)(86,40)(90,40)
\qbezier(90,40)(120,40)(120,20)

\qbezier(20,20)(20, 0)(50, 0)
\qbezier(50, 0)(54, 0)(60, 3)
\qbezier(60, 3)(70, 8)(80, 3)
\qbezier(80, 3)(86, 0)(90, 0)
\qbezier(90, 0)(120, 0)(120,20)

\qbezier(40,20)(50,10)(60,20)
\qbezier(40,20)(50,30)(60,20)

\qbezier(80,20)(90,10)(100,20)
\qbezier(80,20)(90,30)(100,20)

\thinlines

\multiput(0,0)(40,0){3}
{
\qbezier(20,20)(20,25)(30,25)
\qbezier(30,25)(40,25)(40,20)
\qbezier[10](20,20)(20,15)(30,15)
\qbezier[10](30,15)(40,15)(40,20)
}

\multiput(0,0)(40,0){2}
{
\qbezier(32,20)(32,30)(50,30)
\qbezier(50,30)(68,30)(68,20)
\qbezier(32,20)(32,10)(50,10)
\qbezier(50,10)(68,10)(68,20)
}

\end{picture}
\caption{Dehn twists generating the mapping class group of $\GTS$.}
\label{genericlickgen}
\end{figure}
That the induced map is one-to-one is more difficult. 
A proof can be found in \cite{Birman}.
\end{proof}

\begin{prop}
\label{gts_sphere}
The quotient map $\DiffZ \GTS \rightarrow \Diff(S^2,6)$
induces a short exact sequence
\[
0 
\rightarrow \Zed_2 
\rightarrow \pi_0 \DiffZ\GTS 
\rightarrow \pi_0 \Diff(S^2,6)
\rightarrow 0,
\]
where the generator of $\Zed_2$
is mapped to the standard involution of $\GTS$.
\end{prop}

\begin{proof}
Onto is easy:
Each of the five Dehn twists shown in Figure \ref{genericlickgen}
is sent to a half Dehn twist
around a curve separating two puncture points from the rest.
Two examples are shown in Figure \ref{artin1}.
\begin{figure}[ht!]
\centering
\begin{picture}(260,120)(20,-40)

\thicklines

\qbezier(20,20)(20,40)(50,40)
\qbezier(50,40)(54,40)(60,37)
\qbezier(60,37)(70,32)(80,37)
\qbezier(80,37)(86,40)(90,40)
\qbezier(90,40)(120,40)(120,20)

\qbezier(20,20)(20, 0)(50, 0)
\qbezier(50, 0)(54, 0)(60, 3)
\qbezier(60, 3)(70, 8)(80, 3)
\qbezier(80, 3)(86, 0)(90, 0)
\qbezier(90, 0)(120, 0)(120,20)

\multiput(0,0)(40,0){2}
{
\qbezier(40,20)(50,10)(60,20)
\qbezier(40,20)(50,30)(60,20)
}

\thinlines

\qbezier(20,20)(20,25)(30,25)
\qbezier(30,25)(40,25)(40,20)
\qbezier[10](20,20)(20,15)(30,15)
\qbezier[10](30,15)(40,15)(40,20)

\qbezier(72,20)(72,30)(90,30)
\qbezier(90,30)(108,30)(108,20)
\qbezier(72,20)(72,10)(90,10)
\qbezier(90,10)(108,10)(108,20)

\put(130,20){\line(1,0){20}}
\drawline(150,20)(147,23)
\drawline(150,20)(147,17)

\put(220,20){\bigcircle{120}}
\qbezier(160,20)(160,10)(220,10)
\qbezier(220,10)(280,10)(280,20)
\qbezier[20](160,20)(160,30)(220,30)
\qbezier[20](220,30)(280,30)(280,20)
\multiput(190,50)(12,0){6}{\circle*{4}}
\multiput(196,50)(36,0){2}{\bigcircle{24}}
\end{picture}
\caption{Dehn twists mapped to half Dehn twists.}
\label{artin1}
\end{figure}
The definition of a half Dehn twist is as illustrated in Figure \ref{artin2}.
These half Dehn twists are the standard generators
of the mapping class group of the $6$-times punctured sphere.
\begin{figure}[ht!]
\centering
\begin{picture}(160,60)(0,0)

\put(30,30){\bigcircle{60}}
\multiput(20,30)(20,0){2}{\circle*{4}}
\put(30,0){\line(0,1){60}}

\put(70,30){\line(1,0){20}}
\drawline(90,30)(87,33)
\drawline(90,30)(87,27)

\put(130,30){\bigcircle{60}}
\multiput(120,30)(20,0){2}{\circle*{4}}
\qbezier(130,60)(130,55)(120,45)
\qbezier(120,45)(110,35)(110,30)
\qbezier(110,30)(110,20)(120,20)
\qbezier(120,20)(130,20)(130,30)
\qbezier(130,30)(130,40)(140,40)
\qbezier(140,40)(150,40)(150,30)
\qbezier(150,30)(150,25)(140,15)
\qbezier(140,15)(130,5)(130,0)
\end{picture}
\caption{A half Dehn twist}
\label{artin2}
\end{figure}

That the kernel is $\Zed_2$ is 
an elementary exercise in (branched) covering space theory.
\end{proof}

In Section \ref{SPHERE}
we constructed a faithful representation
$\mathcal K_n$ of $\pi_0 \Diff(S^2,n)$.
By the previous two propositions,
$\mathcal K_6$
is a representation of $\pi_0 \Diff\GTS$
whose kernel is equal to $\Zed_2$,
generated by the standard involution.

Let $\mathcal H$ be the representation of $\pi_0 \Diff\GTS$
induced by the action of $\Diff\GTS$ on $H_1\GTS$.
This is called the {\em symplectic representation}.
Under this representation, the standard involution is sent to $-I$.
The direct sum $\mathcal K_6 \oplus \mathcal H$
is therefore a faithful representation of $\pi_0 \Diff\GTS$.
It has dimension $6\binom{5}{2} + 4 = 64$.

\section{Matrices}
\label{MATRICES}

We start off by computing matrices for the representation
$\mathcal L_n'$. Explicit matrices
for $\mathcal L_n$ were worked out both in Krammer and Bigelow's work.
We use the conventions of \cite{Bigelow},
but we correct a sign error which occurs in that paper.
Here, $\sigma_i$ are the half Dehn twist generators of the
mapping class group of a punctured disk,
and $\mathcal L_n(\sigma_i)$
acts on the vector space $V$ with basis
$v_{j,k}$ for $1 \leq j < k \leq n$.
$$
\mathcal L_n(\sigma_i)v_{j,k} =
\left\{
\begin{array}{ll}
v_{j,k} & i\notin \{j-1,j,k-1,k\}, \\
qv_{i,k} + (q^2-q)v_{i,j} + (1-q)v_{j,k} & i=j-1 \\
v_{j+1,k} & i=j\neq k-1, \\
qv_{j,i} + (1-q)v_{j,k} - 
   (q^2-q)tv_{i,k} & i=k-1\neq j,\\
v_{j,k+1} & i=k,\\
-tq^2v_{j,k} & i=j=k-1.
\end{array}
\right.
$$
Using this, we can compute
$\exp \circ \ab(\sigma_i) = t^{-1/d}q^{-n/d}$, with $d=\binom{n}{2}$.
Consequently, $\mathcal L_n'(\sigma_i)= t^{-1/d}q^{-n/d} \mathcal L_n(\sigma_i)$.


The induced representation $\mathcal K_n$ of 
$\mathcal L_{n-1}'$ is now 
straightforward to compute, and we will give 
a block-matrix description of it in terms of $\mathcal L_{n-1}'$.

Reminder: suppose a subgroup $A$ of a group $B$
acts on a vector space $V$.
The {\em induced representation} of $B$
is the module $\Map^A(B,V)$ of $A$-equivariant maps from $B$ to $V$.
The action of $B$ on this module is given by
$b.f := f \circ R_b$,
where $R_b \co B \to B$ is right multiplication by $b$.
Let $\{c_i\}$ be a set of coset representatives of $A$ in $B$,
ie., $B$ is the disjoint union of the cosets $c_i A$.
Then $\Map^A(B,V) = \oplus_i c_i.V$, where our inclusion 
$V \rightarrowtail \Map^A(B,V)$ is given by the $A$-equivariant maps from 
$B$ to $V$ which are zero outside of $A$. The 
direct sum is in the category of abelian groups. 
See \cite[Proposition XVIII.7.2]{Lang} for details.

As coset representatives for $\Stab{p_n}$ in $\pi_0 \Diff(S^2,n)$ we 
will use the maps
$c_1=Id$, 
$c_2=\sigma_{n-1}$, and
$$c_i=(\sigma_{n-i+1} \sigma_{n-i+2} \dots \sigma_{n-2}) \sigma_{n-1}
      (\sigma_{n-i+1} \sigma_{n-i+2} \dots \sigma_{n-2})^{-1}$$
for $i=3,\dots,n$.
Let $\phi_i$ be the permutation of $\{1,\dots, n\}$ such that
$\sigma_i c_j$ is in the coset $c_{\phi_i j} \Stab{p_n}$.
Thus $\phi_i$ is the transposition $(n-i,n-i+1)$.
Then
$$\sigma_i(c_j.v) = c_{\phi_i j}.(c_{\phi_i j}^{-1} \sigma_i c_j v),$$
for any $i=1,\dots,n-2$, $j=1,\dots,n$ and $v \in V$.

Let $\tau = \sigma_1 \sigma_2 \dots \sigma_{n-2} \sigma_{n-2} \dots \sigma_2 \sigma_1$
and let $\nu_j=\sigma_{n-j+1}\sigma_{n-j+2}\dots \sigma_{n-2}$.
Then:
$$
c_{\phi_ij}^{-1}\sigma_ic_{j}=
\left\{
\begin{array}{ll}
 \sigma_i & i\neq n-1, j \neq n+1-i \\
 (\sigma_1 \dots \sigma_{i-1}) \tau^{-1} (\sigma_1 \dots \sigma_{i-1})^{-1} \sigma_i^{-1} & i\neq n-1, j = n+1-i \\
 Id & i=n-1, j=1 \\
 \sigma_{n-2}\tau^{-1} & i=n-1, j=2 \\
 \nu_j\sigma_{n-2}\nu_j^{-1} & i=n-1, j>2
\end{array}
\right.
$$
One can now deduce the matrices $\mathcal K_n(\sigma_i)$.

\section{Remarks}
\label{REMARKS}

Equipped with the knowledge that the mapping class group of an arbitrarily
punctured sphere is linear, Theorem 1 from \cite{Birman} 
allows us to deduce that several subgroups of mapping
class groups are linear.

Let $S$ be a closed $2$-manifold
together with a group $G$ of covering transformations acting on it.
The covering transformations are allowed to have
a finite number of branch points.
Let $n$ be the number of branch points of the covering space $S \to S/G$
and let $\Diff^G S$ be 
the group of fiber-preserving diffeomorphisms of that covering space.
An easy covering space argument shows that
there is an exact sequence of groups
$$G \to \pi_0 \Diff^G S \to \pi_0 \Diff(S/G,n).$$
Suppose there is a faithful representation of $\pi_0 \Diff(S/G,n)$.
Then the above exact sequence gives a representation of $\pi_0 \Diff^G S$
whose kernel is the image of $G$.
If $G$ acts faithfully on $H_1(S)$
then we can obtain a faithful representation of $\pi_0 \Diff^G S$
by taking a direct sum with the symplectic representation.

Suppose $G$ is solvable and fixes each branch point,
and $S$ is not a sphere or a torus.
Then \cite[Theorem 1]{Birman}
states that the map $\Diff^G S \to \Diff S$
induces an injection $\pi_0 \Diff^G S \to \pi_0 \Diff S$.
We claim that $\pi_0 \Diff^G S$ is
the normalizer of $G$ in $\pi_0 \Diff S$.
The proof of this claim uses the fact that any element of $\pi_0 \Diff S$ 
which normalizes the image of $G$ in $\pi_0 \Diff S$
can be lifted to an element of $\Diff S$ which normalizes $G$.
This is proved for the case $G$ is cyclic in \cite[Theorem 3]{Birman}.
The general case follows exactly the same proof
but uses the fact that the Nielsen realization problem
is now solved for all finite groups \cite{Kerckhoff}.

The above line of reasoning can be used
to obtain a faithful representation
of the {\em hyperelliptic mapping class group} of a closed surface $S$.
This is the group of elements of $\pi_0 \Diff S$
which commute with the hyperelliptic involution.
In this case the group $G$ is $\Zed_2$, 
generated by the hyperelliptic involution.
The quotient $S/G$ is a sphere with $2g+2$ branch points.
The generator of $G$ acts as $-I$ on $H_1(S)$.

More generally,
if $S \to S^2$ is a branched covering space
such that the group of covering transformations is solvable
and fixes the branch points
then the normalizer of $G$ in $\Diff(S)$ is linear.
The argument proceeds as previously
except we need to show that $G$ acts faithfully on $H_1(S)$.
This follows from the well-known fact that the Torelli group is torsion-free.
One way to see this is to realize
a torsion element as an isometry of the surface with a suitable hyperbolic
structure \cite{Kerckhoff}. Such a map cannot be trivial on homology
(see, for example \cite[Section V.3]{Farkas}).

Finally, note that if $S$ is a finite-sheeted covering space of $\GTS$
without branch points, with solvable group of covering transformations,
then by the same methods, we obtain a faithful representation of the 
normalizer of the group of covering transformations in $\pi_0 \Diff S$.

\providecommand{\bysame}{\leavevmode\hbox to3em{\hrulefill}\thinspace}

\Addresses
\recd

\end{document}

%% file: agtout.tex
\def\ifplaintex{\expandafter\ifx\csname documentclass\endcsname\relax}

\def\gtp{{\mathsurround=0pt\it $\cal G\mskip-2mu$eometry \&\ 
$\cal T\!\!$opology $\cal P\!$ublications}}  

\def\recd{{\small Received:\qua\receiveddate\ifx\reviseddate\relax
\else\qquad Revised:\qua\reviseddate\fi\par}} 

\def\lognumber#1{\def\thelognumber{#1}}
\def\volumenumber#1{\def\thevolumenumber{#1}}
\def\volumeyear#1{\def\thevolumeyear{#1}}
\def\papernumber#1{\def\thepapernumber{#1}}
\def\pagenumbers#1#2{\def\startpage{#1}\def\finishpage{#2}}
\def\published#1{\def\publishdate{#1}}

\def\received#1{\def\receiveddate{#1}}
\def\revised#1{\def\reviseddate{#1}}
\def\accepted#1{\def\accepteddate{#1}}
\def\asciititle#1{\def\theasciititle{#1}}

\def\asciiaddress#1{\def\theasciiaddress{#1}}
\def\asciiemail#1{\def\theasciiemail{#1}}

\long\def\asciiabstract#1{\long\def\theasciiabstract{#1}}

\let\\\par\let\thelognumber\relax\let\thevolumenumber\relax
\let\thepapernumber\relax\let\thevolumeyear\relax\let\startpage\relax
\let\finishpage\relax\let\publishdate\relax\let\receiveddate\relax
\let\reviseddate\relax\let\accepteddate\relax\let\theasciititle\relax
\let\theasciiauthors\relax\let\theasciiaddress\relax
\let\theasciiabstract\relax

\let\theasciiemail\relax

\ifplaintex
\font\logobig=cmssbx10 scaled 3836
\font\logomed=cmssbx10 scaled 2557
\else
\font\logobig=cmssbx10 scaled 4200
\font\logomed=cmssbx10 scaled 2800
\fi

\long\def\makeagttitle{   
\count0=\startpage
\agt\hfill      
\hbox to 45truept{\vbox to 0pt{\vglue -13truept{\logomed A\kern -.37em{\logobig 
T}\kern -.38em G}\vss}\hss}
\break
{\small Volume \thevolumenumber\ (\thevolumeyear)
\startpage--\finishpage\nl
Published: \publishdate}

\vglue .25truein

{\parskip=0pt\leftskip 0pt plus
1fil\def\\{\par\smallskip}{\Large\bf\thetitle}\par\medskip} \vglue
0.05truein

{\parskip=0pt\leftskip 0pt plus 1fil\def\\{\par}{\sc\theauthors}
\par\medskip}%
 
\vglue 0.03truein 

{\small\leftskip 25truept\rightskip 25truept{\bf Abstract}\stdspace\theabstract

{\bf AMS Classification}\stdspace\theprimaryclass
\ifx\thesecondaryclass\relax\else; \thesecondaryclass\fi\par
{\bf Keywords}\stdspace \thekeywords\par}\vglue 7truept

}   

\ifplaintex
\font\phead=cmsl9 scaled 950
\font\pnum=cmbx10 scaled 913
\font\pfoot=cmsl9 scaled 950
\headline{\vbox to 0pt{\vskip -4.5mm\line{\small\phead\ifnum
\count0=\startpage ISSN 1472-2739 (on-line) 1472-2747 (printed)
\hfill {\pnum\folio}\else\ifodd\count0\def\\{ }%
\ifx\theshorttitle\relax\thetitle\else\theshorttitle\fi\hfill{\pnum\folio}
\else\def\\{ and }{\pnum\folio}\hfill\ifx\theshortauthors\relax\theauthors
\else\theshortauthors\fi\fi\fi}\vss}}
\footline{\vbox to 0pt{\vglue 0mm\line{\small\pfoot\ifnum\count0=\startpage
\copyright\ \gtp\hfill\else
\agt, Volume \thevolumenumber\ (\thevolumeyear)\hfill\fi}\vss}}
\else
\font=cmsl9 scaled 1050
\font\lnum=cmbx10 
\font=cmsl9 scaled 1050
\makeatletter
\def\@oddhead{{\small\lhead\ifnum\count0=\startpage ISSN 1472-2739 
(on-line) 1472-2747 (printed)\hfill {\lnum\number\count0}\else\ifodd\count0
\def\\{ }\ifx\theshorttitle\relax \thetitle \else\theshorttitle\fi\hfill
{\lnum\number\count0}\else\def\\{ and }{\lnum\number\count0}
\hfill\ifx\theshortauthors\relax 
\theauthors\else\theshortauthors\fi\fi\fi}}\def\@evenhead{\@oddhead}
\def\@oddfoot{\small\lfoot\ifnum\count0=\startpage\copyright\ \gtp\hfill\else
\agt, Volume \thevolumenumber\ (\thevolumeyear)\hfill\fi}
\def\@evenfoot{\@oddfoot}
\makeatother
\fi
\let\maketitlepage\makeagttitle
\let\makeshorttitle\maketitlepage
\let\maketitle\maketitlepage

\newwrite\gtoutfile
\long\gdef\makeheadfile{  
{\def\\{, }\def\s{ }
\immediate\openout\gtoutfile head.xxx
\immediate\write\gtoutfile{To: math@arxiv.org}
\immediate\write\gtoutfile{Subject: put OR rep NNNNN:ppppp}
\immediate\write\gtoutfile{--text follows this line--}
\immediate\write\gtoutfile{Proxy-for: \ifx\theasciiauthors\relax
\theauthors\else\theasciiauthors\fi\s<\ifx\theasciiemail\relax\theemail\else\theasciiemail\fi>}
\immediate\write\gtoutfile{\noexpand\\}
\immediate\write\gtoutfile{Authors: \ifx\theasciiauthors\relax
\theauthors\else\theasciiauthors\fi}
{\def\\{ }\immediate\write\gtoutfile{Title: \ifx\theasciititle\relax
\thetitle\else\theasciititle\fi}}
\immediate\write\gtoutfile{Subj-class: GT or SG, GR etc}
\immediate\write\gtoutfile{MSC-class: \theprimaryclass\ifx\thesecondaryclass\relax\else, \thesecondaryclass\fi}
\immediate\write\gtoutfile{Journal-ref: Algebr. Geom. Topol. \thevolumenumber\s
(\thevolumeyear) \startpage-\finishpage}
\immediate\write\gtoutfile{Comments: Published by Algebraic and
Geometric Topology at}
\immediate\write\gtoutfile{\s\s\s  http://www.maths.warwick.ac.uk/agt/AGTVol\thevolumenumber/agt-\thevolumenumber-\thepapernumber.abs.html}
\immediate\write\gtoutfile{\noexpand\\}
\immediate\write\gtoutfile{}
\ifx\theasciiabstract\relax
\immediate\write\gtoutfile{\theabstract}\else
\immediate\write\gtoutfile{\theasciiabstract}\fi
\immediate\write\gtoutfile{}
\immediate\write\gtoutfile{\noexpand\\}
\immediate\write\gtoutfile{}
\immediate\closeout\gtoutfile}}  

\def\maketitlepage{\makeagttitle\makeheadfile}
\let\makeshorttitle\maketitlepage
\let\maketitle\maketitlepage

\def\ifplaintex{\expandafter\ifx\csname documentclass\endcsname\relax}

\def\gtp{{\mathsurround=0pt\it $\cal G\mskip-2mu$eometry \&\ 
$\cal T\!\!$opology $\cal P\!$ublications}}  

\def\recd{{\small Received:\qua\receiveddate\ifx\reviseddate\relax
\else\qquad Revised:\qua\reviseddate\fi\par}} 

\def\lognumber#1{\def\thelognumber{#1}}
\def\volumenumber#1{\def\thevolumenumber{#1}}
\def\volumeyear#1{\def\thevolumeyear{#1}}
\def\papernumber#1{\def\thepapernumber{#1}}
\def\pagenumbers#1#2{\def\startpage{#1}\def\finishpage{#2}}
\def\published#1{\def\publishdate{#1}}

\def\received#1{\def\receiveddate{#1}}
\def\revised#1{\def\reviseddate{#1}}
\def\accepted#1{\def\accepteddate{#1}}
\def\asciititle#1{\def\theasciititle{#1}}

\def\asciiaddress#1{\def\theasciiaddress{#1}}
\def\asciiemail#1{\def\theasciiemail{#1}}

\long\def\asciiabstract#1{\long\def\theasciiabstract{#1}}

\let\\\par\let\thelognumber\relax\let\thevolumenumber\relax
\let\thepapernumber\relax\let\thevolumeyear\relax\let\startpage\relax
\let\finishpage\relax\let\publishdate\relax\let\receiveddate\relax
\let\reviseddate\relax\let\accepteddate\relax\let\theasciititle\relax
\let\theasciiauthors\relax\let\theasciiaddress\relax
\let\theasciiabstract\relax

\let\theasciiemail\relax

\ifplaintex
\font\logobig=cmssbx10 scaled 3836
\font\logomed=cmssbx10 scaled 2557
\else
\font\logobig=cmssbx10 scaled 4200
\font\logomed=cmssbx10 scaled 2800
\fi

\long\def\makeagttitle{   
\count0=\startpage
\agt\hfill      
\hbox to 45truept{\vbox to 0pt{\vglue -13truept{\logomed A\kern -.37em{\logobig 
T}\kern -.38em G}\vss}\hss}
\break
{\small Volume \thevolumenumber\ (\thevolumeyear)
\startpage--\finishpage\nl
Published: \publishdate}

\vglue .25truein

{\parskip=0pt\leftskip 0pt plus
1fil\def\\{\par\smallskip}{\Large\bf\thetitle}\par\medskip} \vglue
0.05truein

{\parskip=0pt\leftskip 0pt plus 1fil\def\\{\par}{\sc\theauthors}
\par\medskip}%
 
\vglue 0.03truein 

{\small\leftskip 25truept\rightskip 25truept{\bf Abstract}\stdspace\theabstract

{\bf AMS Classification}\stdspace\theprimaryclass
\ifx\thesecondaryclass\relax\else; \thesecondaryclass\fi\par
{\bf Keywords}\stdspace \thekeywords\par}\vglue 7truept

}   

\ifplaintex
\font\phead=cmsl9 scaled 950
\font\pnum=cmbx10 scaled 913
\font\pfoot=cmsl9 scaled 950
\headline{\vbox to 0pt{\vskip -4.5mm\line{\small\phead\ifnum
\count0=\startpage ISSN 1472-2739 (on-line) 1472-2747 (printed)
\hfill {\pnum\folio}\else\ifodd\count0\def\\{ }%
\ifx\theshorttitle\relax\thetitle\else\theshorttitle\fi\hfill{\pnum\folio}
\else\def\\{ and }{\pnum\folio}\hfill\ifx\theshortauthors\relax\theauthors
\else\theshortauthors\fi\fi\fi}\vss}}
\footline{\vbox to 0pt{\vglue 0mm\line{\small\pfoot\ifnum\count0=\startpage
\copyright\ \gtp\hfill\else
\agt, Volume \thevolumenumber\ (\thevolumeyear)\hfill\fi}\vss}}
\else
\font=cmsl9 scaled 1050
\font\lnum=cmbx10 
\font=cmsl9 scaled 1050
\makeatletter
\def\@oddhead{{\small\lhead\ifnum\count0=\startpage ISSN 1472-2739 
(on-line) 1472-2747 (printed)\hfill {\lnum\number\count0}\else\ifodd\count0
\def\\{ }\ifx\theshorttitle\relax \thetitle \else\theshorttitle\fi\hfill
{\lnum\number\count0}\else\def\\{ and }{\lnum\number\count0}
\hfill\ifx\theshortauthors\relax 
\theauthors\else\theshortauthors\fi\fi\fi}}\def\@evenhead{\@oddhead}
\def\@oddfoot{\small\lfoot\ifnum\count0=\startpage\copyright\ \gtp\hfill\else
\agt, Volume \thevolumenumber\ (\thevolumeyear)\hfill\fi}
\def\@evenfoot{\@oddfoot}
\makeatother
\fi
\let\maketitlepage\makeagttitle
\let\makeshorttitle\maketitlepage
\let\maketitle\maketitlepage

\newwrite\gtoutfile
\long\gdef\makeheadfile{  
{\def\\{, }\def\s{ }
\immediate\openout\gtoutfile head.xxx
\immediate\write\gtoutfile{To: math@arxiv.org}
\immediate\write\gtoutfile{Subject: put OR rep NNNNN:ppppp}
\immediate\write\gtoutfile{--text follows this line--}
\immediate\write\gtoutfile{Proxy-for: \ifx\theasciiauthors\relax
\theauthors\else\theasciiauthors\fi\s<\ifx\theasciiemail\relax\theemail\else\theasciiemail\fi>}
\immediate\write\gtoutfile{\noexpand\\}
\immediate\write\gtoutfile{Authors: \ifx\theasciiauthors\relax
\theauthors\else\theasciiauthors\fi}
{\def\\{ }\immediate\write\gtoutfile{Title: \ifx\theasciititle\relax
\thetitle\else\theasciititle\fi}}
\immediate\write\gtoutfile{Subj-class: GT or SG, GR etc}
\immediate\write\gtoutfile{MSC-class: \theprimaryclass\ifx\thesecondaryclass\relax\else, \thesecondaryclass\fi}
\immediate\write\gtoutfile{Journal-ref: Algebr. Geom. Topol. \thevolumenumber\s
(\thevolumeyear) \startpage-\finishpage}
\immediate\write\gtoutfile{Comments: Published by Algebraic and
Geometric Topology at}
\immediate\write\gtoutfile{\s\s\s  http://www.maths.warwick.ac.uk/agt/AGTVol\thevolumenumber/agt-\thevolumenumber-\thepapernumber.abs.html}
\immediate\write\gtoutfile{\noexpand\\}
\immediate\write\gtoutfile{}
\ifx\theasciiabstract\relax
\immediate\write\gtoutfile{\theabstract}\else
\immediate\write\gtoutfile{\theasciiabstract}\fi
\immediate\write\gtoutfile{}
\immediate\write\gtoutfile{\noexpand\\}
\immediate\write\gtoutfile{}
\immediate\closeout\gtoutfile}}  

\def\maketitlepage{\makeagttitle\makeheadfile}
\let\makeshorttitle\maketitlepage
\let\maketitle\maketitlepage